\documentclass[preprint,3p]{elsarticle}
\usepackage{amsmath,amsfonts}
\usepackage[matrix,arrow,curve]{xy}
\usepackage{graphicx}
\usepackage{pstricks, pst-plot, pst-node}
\usepackage{subfigure}
\usepackage{enumerate}

\def\subheading#1{\medskip\noindent{\boldmath\textbf{#1.} }\ignorespaces}

\def\F{{\bf F}}

\newcommand{\R}{\mathbb{R}}

\newcommand{\Ocal}{\mathcal{O}}

\newcommand{\divi}{\mathrm{div}}

\newcommand{\ie}{i.\,e.\ }

\def\proj{{\mathbb{P}}}

\def\url#1{#1}

\newtheorem{theorem}{Theorem}

\newtheorem{example}[theorem]{Example}
\newtheorem{remark}[theorem]{Remark}

\newproof{proof}{Proof}

\def\EE#1#2{{\rm E}_{#1,#2}}

\def\ssquarings{{\bf s}} %squarings over small field
\def\smults{{\bf m}} %multiplications over small field
\def\smultsbya{{\bf m_a}} %multiplications by parameter a
\def\smultsbyd{{\bf m_d}} %multiplications by parameter d

\def\squarings{{\bf S}}
\def\mults{{\bf M}}

\def\Ocal{\mathcal{O}}
\def\Jcal{\mathcal{J}}
\def\Pcal{\mathcal{P}}
\def\Ecal{\mathcal{E}}
\def\ie{i.\,e.\ }

\begin{document}
\title{Faster Computation of the Tate Pairing}

\ifodd0
\author{ }
\institute{ }
\fi

\ifodd1
\author[iml]{Christophe Ar\`ene}
\ead{arene@iml.univ-mrs.fr}

\author[eindhoven]{Tanja Lange\corref{cor1}}
\ead{tanja@hyperelliptic.org}
\cortext[cor1]{Corresponding author}
\ead[url]{hyperelliptic.org/tanja}

\author[eindhoven,microsoft]{Michael Naehrig}
\ead{michael@cryptojedi.org}
\ead[url]{cryptojedi.org/michael}

\author[iml]{Christophe Ritzenthaler}
\ead{ritzenth@iml.univ-mrs.fr}

\address[iml]{
Institut de Math\'ematiques de Luminy\\
163, avenue de Luminy, Case 907\\
13288 Marseille CEDEX 09\\
France}

\address[eindhoven]{
Department of Mathematics and Computer Science\\
Technische Universiteit Eindhoven\\
P.O. Box 513, 5600 MB Eindhoven\\
Netherlands}

\address[microsoft]{
Microsoft Research\\
One Microsoft Way\\
Redmond, 98052 WA\\
USA}

\ifodd1
\catcode`@=11
\def\@thefnmark{\ }
\@footnotetext{
This work has been supported in part by the European Commission
through the ICT Programme under Contract ICT--2007--216646 ECRYPT II,
and in part by grant MTM2006-11391 from the Spanish MEC. The first author is
beneficiary of a Ph.D. grant from the AXA Research Fund.}
\catcode`@=12
\fi

\begin{abstract}
  This paper proposes new explicit formulas for the doubling and
  addition steps in Miller's algorithm to compute the Tate pairing on elliptic
  curves in Weierstrass and in Edwards form. 
  For Edwards curves the formulas come from a new way of seeing 
the arithmetic. 
  We state the first geometric interpretation of the group law
  on Edwards curves by presenting the functions which arise in
  addition and doubling. The Tate pairing on Edwards curves can be computed by
  using these functions in Miller's algorithm.
  
  Computing the sum of two points or the double of a point and the coefficients
  of the corresponding functions is faster with our formulas than with all previously proposed
  formulas for pairings on Edwards curves. They are even competitive with
  all published formulas for pairing computation on Weierstrass curves.
%% Check loop shortening before submitting:
%This idea can be combined with loop shortening, as long as the
%  pairing is based on the Tate pairing and the curve
%  has order divisible by $4$.
  We also improve the formulas for Tate pairing computation on Weierstrass curves in
  Jacobian coordinates.
Finally, we present several examples of pairing-friendly Edwards curves.
\end{abstract}


\maketitle

\section{Introduction}\label{intro}
Since their introduction to cryptography by Bernstein and
Lange~\cite{2007/BernsteinLange}, Edwards curves have received a lot
of attention due to the fact that their group law can be computed very
efficiently. The group law in
affine form was introduced by Edwards in \cite{2007/edwards} along
with a description of the curve and several proofs of correctness.
Remarkably none of the proofs provided a geometric interpretation
while addition on   Weierstrass curves
is usually explained via the chord-and-tangent method. 

Cryptographic applications in discrete-logarithm-based systems such as
Diffie-Hellman key exchange or digital signatures require efficient
computation of scalar multiples and thus have benefited from the speedup in
addition and doubling. The situation is significantly different in
pairing-based cryptography where Miller's algorithm needs 
a function with divisor $(P)+(Q)-(P+Q)-(\Ocal)$ for two 
input points $P$ and $Q$, their sum $P+Q$, and neutral element $\Ocal$. For curves in Weierstrass
form these functions are readily given by the line functions in the
usual addition and doubling. Edwards curves have degree $4$ and thus
any line passes through $4$ curve points instead of $3$. This led many to
conclude that Edwards curves provide no benefit to pairings and are
doomed to be slower than the Weierstrass counterparts.

So far two papers have attempted to compute pairings efficiently 
on Edwards
curves: Das and Sarkar~\cite{2008/das} use the birational equivalence
to Weierstrass curves to map the points on the Edwards curve to a
Weierstrass curve on which the usual line functions are then
evaluated. This approach comes at a huge performance penalty as these
implicit pairing formulas need many field operations to evaluate
them. Das and Sarkar then focus on supersingular curves with embedding
degree $k=2$ and develop explicit formulas for that case.
% but they argue that the birational equivalence would take the affine
% point on the Edwards curve to a projective point on the Weierstrass
% curve so that either projective addition formulas or an inversion need
% to be used which makes 

Ionica and Joux~\cite{2008/ionica} use a different map to a curve of
degree $3$ and compute the $4$-th power of the Tate pairing. The
latter poses no problem for usage in protocols as long as all participating
parties
perform the same type of pairing computation. Their results are
significantly faster than Das and Sarkar's but they are still much slower
than pairings on Weierstrass curves.

In this paper we close several important gaps:
\begin{itemize}
\item We provide a geometric interpretation of the addition
  law for twisted Edwards curves.
\item We study additions, doublings, and all the special cases that
  appear as part of the geometric addition law for twisted Edwards
  curves.
\item We use the geometric interpretation of the group law to show how to compute the Tate
  pairing on twisted Edwards curves.
\item We give examples of ordinary pairing-friendly Edwards
  curves at several security levels. The curves have embedding degrees between
  $6$ and $22$.
\end{itemize}

Beyond that, we develop explicit formulas for computing the Tate pairing on
Edwards curves that
\begin{itemize}
\item solidly beat the results by Das and Sarkar~\cite{2008/das} and
  Ionica and Joux~\cite{2008/ionica};
%\item are faster than the fastest formulas for computations on
%  Weierstrass curves in Jacobian coordinates 
%  (\cite{2003/izu} and \cite{2008/ionica});
%\item are faster than the fastest formulas for computations on
%  Weierstrass curves with $a_4=-3$ in Jacobian coordinates
%  (\cite{2004/chatterjee} and \cite{2003/izu});
\item are as fast as the fastest previously published formulas for the doubling
  step on Weierstrass curves, namely curves with $a_4=0$ (e.g.
  Barreto-Naehrig curves) in Jacobian coordinates, and faster than
  other Weierstrass curves;
\item need the same number of field operations as the best published
  formulas for mixed addition in Jacobian coordinates; and  
\item have minimal performance penalty for non-affine base points.
\end{itemize}

In particular, for even embedding degree $k$ the doubling step on an
Edwards curve takes $1\mults +1\squarings+(k+6)\smults+5\ssquarings$,
where $\smults$ and $\ssquarings$ denote the costs of multiplication
and squaring in the base field while $\mults$ and $\squarings$ denote
the costs of multiplication and squaring in the extension field of
degree $k$. A mixed addition step takes $1\mults+(k+12)\smults$ and an
addition step takes $1\mults+(k+14)\smults$.
Our method for pairing computation on Edwards curves can be used for all curves
that can be represented in Edwards form over the base field.

We also improve the addition and doubling steps on Weierstrass curves given by
an equation $y^2 = x^3 + a_4x + a_6$.
We present the first explicit formulas for
full addition steps on Weierstrass curves. The new formulas need $1\mults
+1\squarings+(k+6)\smults+5\ssquarings$ for a doubling step on curves with
coefficient $a_4 = -3$. On such curves a mixed addition step costs
$1\mults+(k+6)\smults + 6\ssquarings$ and an addition step costs $1\mults+(k +
9)\smults + 6\ssquarings$. On curves with $a_4 = 0$, the formulas take $1\mults
+1\squarings+(k+3)\smults+8\ssquarings$ for a doubling step,
$1\mults+(k+6)\smults + 6\ssquarings$ for a mixed addition, and
$1\mults+(k+9)\smults + 6\ssquarings$ for an addition step.

% The following is not clear. For an s/m ratio of 0.8, 9m+6s is faster than
% 14m. I would rephrase this. (Michael, 2009-09-08)
Our new formulas for Weierstrass curves are the fastest when using affine
base points (except in the case $a_4=0$, $a_6 = b^2$). For projective base
points -- a common case in pairing-based protocols -- it is better to use
Edwards curves.

%This paper does not consider other pairings such as the ate pairing;
%such pairings are particularly interesting for curves with twists of
%degree $4$ or $6$ and our example curves do not fit with that.
%The geometric interpretation and the explicit formulas work
%over any field.
%% Not 100% sure about that, check before submitting:
%It can also be used with loop shortening to compute the
%eta or twisted ate pairing. This is particularly interesting for curves with
%twists of degree $4$ and $6$, where the twisted ate pairing usually leads to
%a shorter Miller loop compared to the Tate pairing.

%{\bf Notation.} In the rest of the paper, $K$ will be a field of characteristic different from $2$.

\section{Background on Pairings}\label{pair}
Let $q$ be a prime power not divisible by $2$ and let $E/\F_q$ be an elliptic curve over $\F_q$ with neutral element
denoted by $\Ocal$. Let $n\mid \#E(\F_q)$ be a prime divisor of the
group order and let $E$ have embedding degree $k>1$ with respect to $n$, \ie $k$
is the smallest integer such that $n \mid q^k-1$.
%For simplicity and speed we assume that $k>1$.

Let $P \in E(\F_q)[n]$ and let $f_P \in \F_q(E)$ be  such that
$\divi(f_P)=n(P)-n(\Ocal)$. Let
$\mu_n\subset \F^*_{q^k}$ denote the group of $n$-th roots of unity.
The reduced Tate pairing is given by
% \begin{eqnarray*}
% T_n &: E(\F_p)[n] \times E(\F_{p^k})/nE(\F_{p^k}) & \to \mu_n \\
%  & (P,Q)& \mapsto f_P(Q)^{(p^k-1)/n}.
%  \end{eqnarray*}
$$T_n: E(\F_q)[n] \times E(\F_{q^k})/nE(\F_{q^k})\to \mu_n;
\ (P,Q)\mapsto f_P(Q)^{(q^k-1)/n}.$$
% put something on the definition of ate pairing here.
Miller~\cite{2004/mill} suggested to compute pairings in an iterative
manner. Let $n=(n_{l-1},\dots, n_1,n_0)_2$ be the binary
representation of $n$, where $n_{l-1} = 1$. Let $g_{R,S}\in \F_q(E)$ be the function
arising in the addition of two points $R$ and $S$ on $E$, i.e. $g_{R,S}$ is a function with
$\divi(g_{R,S})=(R)+(S)-(R+S)-(\Ocal)$, where 
$\Ocal$ denotes the neutral element in the group of points,
$R+S$ denotes the sum of
$R$ and $S$ on $E$, and  additions of the form $(R)+(S)$ denote formal
additions in the divisor group. Miller's algorithm starts with
$R=P, f=1$ and computes
\begin{enumerate}
\item for $i=l-2$ to $0$ do
\begin{enumerate}
\item $f\gets f^2 \cdot g_{R,R}(Q),\ R\gets [2]R$, \hfill //doubling step
\item if $n_i=1$ then $f\gets f \cdot g_{R,P}(Q),\ R\gets R+P$.\hfill //addition step
\end{enumerate}
\item $f\gets f^{(q^k-1)/n}$.
\end{enumerate}

Note that pairings can be combined with windowing methods by replacing
the computation in step~(b) by
$$f\gets f \cdot f_{c,P}(Q) \cdot g_{R,[c]P}(Q),\ R\gets R+[c]P,$$
where the current window in the binary representation of $n$ corresponds to the value $c$.
The Miller function $f_{c,P}$ is defined via $\divi(f_{c,P}) = c(P) -
([c]P) - (c-1)(\Ocal)$.
But windowing methods are rarely
used because of the extra costs of $1\mults$ for updating the variable $f$.

\section{Formulas for Pairings on Weierstrass curves}\label{appendix}
An elliptic curve over $\F_q$ in short Weierstrass form is given by an
equation of the form $y^2=x^3+a_4x+a_6$ with $a_4, a_6\in \F_q$.  In
this section we present new formulas for the addition and doubling
step in Miller's algorithm that are faster than previous ones.
Furthermore, we also cover the case of a non-affine base point.

%  Usual way of speeding up computation in case
% $k$ even is to start with a point $(x_0,y_0)$ over $\F_{p^{k/2}}$ and
% then twist it.
 
%T: Moved the following to the next section.
% On twisted Edwards curves $\EE{a}{d}$, twists affect the
% $x$-coordinate. Let $\F_{p^k}$ have basis $\{1,\alpha\}$ over
% $\F_{p^{k/2}}$ with $\alpha^2=\delta \in \F_{p^{k/2}}$ and let $Q'=(x_0,y_0)\in
% \EE{a\delta}{d\delta}(\F_{p^{k/2}})$.  Twisting $Q'$ with $\alpha$
% ensures that the second argument of the pairing is on
% $\EE{a}{d}(\F_{p^k})$ (and no smaller field) and is of the form
% $Q=(x_0\alpha,y_0)$, where $x_0, y_0 \in \F_{p^{k/2}}$.

The fastest formulas for doublings on Weierstrass curves are given in
Jacobian coordinates (cf. the EFD~\cite{2007/EFD}). A point is
represented as $(X_1:Y_1:Z_1)$ which for $Z_1\neq 0$ corresponds to
the affine point $(x_1,y_1)$ with $x_1=X_1/Z_1^2$ and $y_1=Y_1/Z_1^3$.
To obtain the full speed of pairings on Weierstrass curves it is
useful to represent a point by $(X_1:Y_1:Z_1:T_1)$ with $T_1=Z_1^2$.
This allows one $\ssquarings - \smults$ tradeoff in the addition step
compared with the usual representation $(X_1:Y_1:Z_1)$. If the
intermediate storage is an issue or if $\ssquarings$ is not much
smaller than $\smults$, $T_1$ should not be cached. We present the
formulas including $T_1$ below; %to have the best operation count for
%Weierstrass curves in the comparison
the modifications to omit $T_1$
are trivial.

For $S\in\{R,P\}$, the function $g_{R,S}$ for Weierstrass curves is given as the fraction of the
usual line functions by
$$
g_{R,S}(X:Y:Z)%=\frac{l_{R,P}}{l_{R+P}}
=\frac{(YZ^3_0-Y_0Z^3)-\lambda(XZ^2_0-X_0Z^2)ZZ_0}
{(X-cZ^2)Z},
$$
where $\lambda$ is the slope of the line through $R$ and $S$ (with
multiplicities), $(X_0:Y_0:Z_0)$ is a point on the line, and $c$ is
the $x$-coordinate of $R+S$. When one computes the Tate pairing, the point
$(X_0:Y_0:Z_0)$ and the constants $\lambda$ and $c$ are defined over
the base field $\F_q$.  The function is evaluated at a point
$Q=(X_Q:Y_Q:Z_Q)$ defined over $\F_{q^k}$.

We assume that $k$ is even. This allows us to use several improvements and
speedups that are presented in \cite{2002/barreto} and \cite{2004/barreto}.
As usual, let the field extension $\F_{q^k}$ be constructed via a quadratic
subfield as $\F_{q^k}= \F_{q^{k/2}}(\alpha)$, with $\alpha^2=\delta$
for a non-square $\delta\in\F_{q^{k/2}}$;
and let $Q$ be chosen to be of the form $Q=(x_Q:y_Q\alpha:1)$ with $x_Q,
y_Q\in \F_{q^{k/2}}$. The latter is 
enforced by choosing a point $Q'$ on a quadratic twist of $E$ over
$\F_{q^{k/2}}$ and defining $Q$ as the image of $Q'$ under the twist
isomorphism. The denominator of $g_{R,S}(Q)$ is given by $x_Q -
c$ which is defined over the subfield $\F_{q^{k/2}}$. Thus
only the numerator needs to be considered as all
multiplicative contributions from proper subfields of $\F_{q^k}$ are mapped to
$1$ by the final exponentiation and can be
discarded. 
Furthermore, for addition and doubling in
Jacobian coordinates we can write $\lambda = L_1/Z_3$, where $Z_3$ is the $z$-coordinate of $R+S$ and $L_1$
depends on $R$ and $S$. Since $Z_3$ is defined over $\F_q$, we can instead
compute
$
Z_3(y_QZ^3_0\alpha - Y_0)-L_1(x_QZ^2_0-X_0)Z_0
$
giving $g_{R,S}$ up to factors from subfields  of $\F_{q^k}$.

%For Weierstrass curves and even $k$, several improvements and speedups
%are presented in \cite{2002/barreto} and \cite{2004/barreto}. In
%particular it is common to eliminate all denominators in $g_{R,S}$ by choosing the
%second point $Q$ such that its $x$-coordinate is in a subfield of
%$\F_{q^k}$. The functions $g_{R,S}$ are defined over $\F_q$ and their
%denominators are functions in $x$ only. Writing
%$g_{R,S}(Q)=h_{R,S}(x_Q,y_Q)/s_{R,S}(x_Q)$ with polynomial functions
%$h_{R,S}$ and $s_{R,S}$, one sees that the complete contribution of all
%terms $s_{R,S}(x_Q)$ will be mapped to $1$ by the final exponentiation if
%$x_Q$ is in a proper subfield of $\F_{q^k}$.  The latter is usually
%enforced by choosing a point $Q'$ on a quadratic twist of $E$ over
%$\F_{q^{k/2}}$ and defining $Q$ as the image of $Q'$ under the twist.
%Note that quadratic twists on Weierstrass curves change only the $y$-coordinate
%but leave the $x$-coordinate invariant.

\subsection{Addition steps}
In Miller's algorithm, all additions involve the base point as one
input point so, when computing the line function, $(X_0:Y_0:Z_0)$ can be
chosen as the base point $P$ and all values depending solely on $P$
and $Q$ can be precomputed at the beginning of the computation. For
additions, $P$ is always stated as the second summand, \ie 
$P=(X_2:Y_2:Z_2)$.  

To enable an $\smults-\ssquarings$
tradeoff we compute $2g_{R,P}(Q)$; this does not change the
result of the computation since $2\in \F_q$. % and $g_{R,P}(Q)$ enters multiplicatively.  
Multiplications with $x_Q$ and $y_Q$ cost
$(k/2)\smults$ each; for $k>2$ it is thus useful to rewrite the line function as
$$
l=Z_3\cdot 2y_QZ^3_2\alpha - 2Z_3\cdot Y_2 - L_1\cdot(2(x_QZ^2_2-X_2)Z_2),
$$
needing $(k+1)\smults$ for precomputed $y'_Q=2y_QZ^3_2\alpha$ and
$x'_Q=2(x_QZ^2_2-X_2)Z_2$. Additionally $1\mults$ is needed to update
the variable $f$ in Miller's algorithm.

\subheading{Full addition} We use Bernstein and Lange's
formulas (``add-2007-bl'') from the EFD~\cite{2007/EFD}. 
%This means that all additions involve the base point $P$
%which is fixed throughout the computation, and
We can
cache all values depending solely on $P$. In particular we precompute
(or cache after the first addition or doubling) $R_2=Y_2^2$ and
$S_2=T_2\cdot Z_2$.
%Independent of the value of $a_4$, all doubling formulas compute
%$Y_2^2$. This means that $R_2=Y_2^2$ can be cached since Miller's
%algorithm starts by computing $2P$.  
The numerator of $\lambda$ is $L_1=D-C$.
%
%\begin{small}
\begin{eqnarray*}
A &=& X_1\cdot T_2;\ B = X_2\cdot T_1;\ C = 2Y_1\cdot S_2;\ D = 
((Y_2+ Z_1)^2-R_2-T_1)\cdot T_1;\\
H &=& B-A;\
     I = (2H)^2;\ 
      J = H\cdot I;\ 
      L_1 = D-C;\
      V = A\cdot I;\\
      X_3 &=& L_1^2-J-2V;\
      Y_3 = L_1\cdot (V-X_3)-2C\cdot J;\
      Z_3 = ((Z_1+Z_2)^2-T_1-T_2)\cdot H;\\
      T_3&=&Z_3^2;\ l = Z_3\cdot y'_Q -(Y_2+ Z_3)^2+R_2+T_3-L_1\cdot x'_Q.
\end{eqnarray*}
%\end{small}
The formulas need $1\mults+(k+9)\smults +6\ssquarings$ to compute the
addition step. To our knowledge this is the first set of formulas for
full (non-mixed) addition. If $\smults$ is not significantly more
expensive than $\ssquarings$, some computations should be performed
differently. In particular, $R_2$ needs not be stored, $D$ is computed
as $D=2Y_2\cdot Z_1\cdot T_1$, $l$ contains the term $-2Y_2\cdot Z_3$
instead of $-(Y_2+ Z_3)^2+R_2+T_3$, and the computation of $Z_3$ can
save some field additions.

If the values $T_1, R_2, S_2, T_2, x'_Q$, and $y'_Q$ cannot be stored,
different optimizations are needed; in particular the line function is
computed as
$$
l=((Z_3 \cdot Z_2)\cdot Z^2_2)\cdot y_Q\alpha - Y_2\cdot Z_3 -
(L_1\cdot Z_2)\cdot Z^2_2\cdot x_Q +X_2\cdot (L_1\cdot Z_2)
$$
and the computation costs end up as $1\mults+(k+17)\smults +6\ssquarings$.
% , where the
% $11\smults +5\ssquarings$.
%  come from the costs of full addition and 
% $(k+6)\smults +1\ssquarings$ account for the computation of $l$.

\subheading{Mixed addition} Mixed addition means that the second input
point is in affine representation. Mixed additions occur in scalar multiplication
 if the base point $P$ is given 
%in affine coordinates 
as
$(x_2:y_2:1)$.

We now state the mixed addition formulas based on
Bernstein and Lange's formulas (``add-2007-bl'') from the
EFD~\cite{2007/EFD}. Mixed additions are the usual case studied for
pairings and the evaluation of the line function in $(k+1)\smults$ is standard. However,
most implementations miss the $\ssquarings-\smults$ tradeoff in the
main mixed addition formulas and do not compute the $T$-coordinate.
\begin{eqnarray*}
B &=&  x_2\cdot T_1;\ D = ((y_2+ Z_1)^2-R_2-T_1)\cdot T_1;\ H = B-X_1;\
     I = H^2;\ 
     E=4I;\
      J = H\cdot E;\\
      L_1 &=& (D-2Y_1);\
      V = X_1\cdot E;\ 
      X_3 = L_1^2-J-2V;\
      Y_3 = r\cdot (V-X_3)-2Y_1\cdot J;\\ 
      Z_3 &=& (Z_1+ H)^2-T_1-I;
      T_3=Z_3^2;\
      l = 2Z_3\cdot y_Q\alpha - (y_2+ Z_3)^2+R_2+T_3 -2L_1\cdot (x_Q -x_2).
\end{eqnarray*}
The formulas need $1\mults+(k+6)\smults +6\ssquarings$ to compute the mixed
addition step. % If operations in $\F_{p^k}$ are not too cumbersome and
%$(y_Q\alpha)^2=y_Q^2\delta$ has been precomputed the line function can
%be computed as $ l = (Z_3+ y_Q\alpha)^2-T_3-y_Q^2\delta - y_2\cdot Z_3
%-r\cdot (x_Q -x_2)$ which leads to a total complexity of 
%$1\mults+(k/2+6)\smults +(k/2+6)\ssquarings$.

\penalty500000

\subsection{Doubling steps}
The main differences between the addition and the doubling formulas
are that the doubling formulas depend on the curve coefficients and that the
%line function must be computed with the input to the doubling function
point $(X_0:Y_0:Z_0)$ appearing in the definition of $g_{R,S}$ is
$(X_1:Y_1:Z_1)$, which is changing at every step. So in
particular $Z_0 \ne 1$ and no precomputations (like $x'_Q$ or $y'_Q$ in
the addition step)
can be done.

For arbitrary $a_4$ the equation of the slope is
$\lambda=(3X_1^2+a_4Z_1^4)/(2Y_1Z_1)=(3X_1^2+a_4Z_1^4)/Z_3$.  Thus
$Z_3$ is divisible by $Z_1$ and we can replace $l$ by $l'=l/Z_1$ which
% it's a bit weird to speak of divisibility with finite field elements
will give the same result for the pairing computation.
The value of 
$$
l'=(Z_3 \cdot Z^2_1)\cdot y_Q\alpha - 2Y_1^2 - L_1\cdot Z^2_1\cdot x_Q
+X_1\cdot L_1
$$
can be computed in at most $(k+3)\smults +1\ssquarings$ for arbitrary $a_4$ and with slightly less
operations otherwise.

The formulas by Ionica and Joux \cite{2008/ionica} take into account the doubling formulas
from the EFD for general Weierstrass curves in Jacobian coordinates. 
We thus present new formulas for the more special curves with $a_4=-3$ and
$a_4=0$. 

% Note that instead of computing $T_3=Z_3^2$ at the end of doubling,
% $T_1=Z_1^2$ could be computed at the beginning. However, additions
% benefit from this extra cached value by an $\ssquarings-\smults$
% tradeoff.

\subheading{Doubling on curves with $a_4=-3$}
The fastest doubling formulas are due to Bernstein (see
\cite{2007/EFD} ``dbl-2001-b'') and need $3\smults +5\ssquarings$ for
the doubling.
\begin{eqnarray*}
     A &=& Y_1^2;\
      B = X_1\cdot A;\ 
      C = 3 (X_1-T_1)\cdot (X_1+T_1);\\
      X_3 &=& C^2-8B;\ 
      Z_3 = (Y_1+Z_1)^2-A-T_1;\ 
      Y_3 = C\cdot (4B-X_3)-8A^2;\\
      l&=&(Z_3 \cdot T_1)\cdot y_Q\alpha - 2A - C\cdot T_1\cdot x_Q
+X_1\cdot C;\ 
      T_3=Z_3^2.
\end{eqnarray*}
The complete doubling step thus takes $1\mults + 1\squarings
+(k+6)\smults +5\ssquarings$. Note that $L_1 = C$.

\penalty-5000

\subheading{Doubling on curves with $a_4=0$}
The following formulas compute a doubling in $1\smults
+7\ssquarings$. Note that without $T_1$ and computing $Z_3 = 2Y_1\cdot
Z_1$ a doubling can be computed in $2\smults +5\ssquarings$ which is
always faster (see \cite{2007/EFD}) but the line functions make use of
$Z_1^2$. Note further that here $L_1=E=3X_1^2$ is particularly simple.
\begin{eqnarray*}
A&=&X_1^2;\ B=Y_1^2;\ C=B^2;\ D = 2((X_1 + B)^2 -A - C);\ E = 3 A;\ G = E^2;\\
X_3&=&G - 2 D;\ Y_3 = E \cdot (D - X_3) - 8 C;\  Z_3 = (Y_1 + Z_1)^2 -B - T_1;\\
      l&=&2(Z_3 \cdot T_1)\cdot y_Q\alpha - 4B - 2E\cdot T_1\cdot x_Q
+(X_1+E)^2-A-G;\ T_3=Z_3^2.\
\end{eqnarray*}
The complete doubling step thus takes $1\mults+1\squarings
+(k+3)\smults +8\ssquarings$.

%\section{Twisted Edwards Curves}\label{edw}
\section{Geometric
interpretation of the group law on twisted Edwards curves}\label{edw}
In this section $K$ denotes a field of characteristic different from $2$.
%They found an addition law and showed that any elliptic curve with
%$\#E(\F_q)$ divisible by $4$ is isogenous to a twisted Edwards curve.
%Let $K$ be a field of characteristic ${\rm char} \neq 2$.
A \emph{twisted Edwards curve} over $K$ is a curve given by an affine equation of the
form
$\EE{a}{d}: ax^2+y^2=1+dx^2y^2$
for $a,d\in K^*$ and $a\neq d$.   
Twisted Edwards curves were introduced by Bernstein et al.  in
\cite{2008/BernsteinBirkner} as a generalization of 
Edwards curves \cite{2007/BernsteinLange} which are included as 
$\EE{1}{d}$. 
%The addition law on $\EE{a}{d}$ found a lot of attention in scalar 
%multiplication.
An addition law on points of the curve $\EE{a}{d}$ is given by
 $$(x_1, y_1)+ (x_2, y_2) =
 \left( \frac{x_1 y_2 + y_1 x_2}{1 + d x_1 x_2 y_1 y_2}, \frac{y_1 y_2
     - a x_1 x_2}{1 - d x_1 x_2 y_1 y_2} \right).$$
The neutral element is $\Ocal=(0,1)$, and the negative of $(x_1,y_1)$
is $(-x_1,y_1)$. The point $\Ocal'=(0,-1)$ has order $2$.  The points
at infinity $\Omega_1=(1:0:0)$ and $\Omega_2=(0:1:0)$ are singular and
blow up to two points each.

Edwards curves received a lot of attention because the above addition can be
computed very efficiently, resulting in highly efficient algorithms to carry out
scalar multiplication, a basic tool for many cryptographic protocols.

The name twisted Edwards curves comes from the fact that the set of
twisted Edwards curves is invariant under quadratic twists while a
quadratic twist of an Edwards curve is not necessarily an Edwards
curve. In particular, let $\delta \in K \setminus K^2$
% be a non-square
and let $\alpha^2=\delta$ for some $\alpha$  in a quadratic
extension $K_2$ of $K$. The map $\epsilon : (x,y)\mapsto (\alpha x,y)$
defines a $K_2$-isomorphism between the twisted Edwards curves 
$\EE{a\delta}{d\delta}$
and  
$\EE{a}{d}$.  Hence, the map $\epsilon$ is the prototype of a quadratic twist. Note that  twists change the $x$-coordinate unlike on Weierstrass
curves where they affect the $y$-coordinate.

%\subsection{Geometric Interpretation of the Group Law}\label{geom}
We now study the intersection of $\EE{a}{d}$ with certain plane curves 
 and explain the Edwards addition law in terms of the
divisor class arithmetic. We remind the reader that the divisor class
group is defined as the group of degree-$0$ divisors modulo the group of
principal divisors in the function field of the curve, i.e. two
divisors are \emph{equivalent} if they differ by a principal divisor. For
background reading on curves and Jacobians, we refer to
\cite{2005/frey-ehcc4} and \cite{1986/silverman}.\\
%We first consider projective
%lines in $\proj^2$. A general line is of the form $L: c_XX + c_YY + c_ZZ = 0$ 
%\begin{equation}\label{eqn:projline}
%\end{equation}
%where $(c_X:c_Y:c_Z) \in \proj^2$. A line is uniquely determined by
%two of its points when they are distinct. 
Let $\proj^2(K)$ be the two-dimensional projective space over $K$, and let $P = (X_0: Y_0: Z_0)
\in \proj^2(K)$ with $Z_0 \neq 0$. Let $L_{1,P}$ be the line through $P$ and
$\Omega_1$, \ie $L_{1,P}$
is defined by $Z_0Y - Y_0Z = 0$; and let $L_{2,P}$ be the line through $P$ and
$\Omega_2$, \ie  $L_{2,P}$ is defined by $Z_0X - X_0Z = 0$.

% We first consider lines
% which pass through one of the points at infinity and an affine point $P$. Note
% that the line through $\Omega_1$ and $\Omega_2$ is the line at infinity
% $L_{\infty}: Z = 0$.

% \begin{lemma}\label{lem:edlines}
% Let $P = (X_0: Y_0: Z_0) \in \proj^2(K)$ be an affine point, \ie $Z_0 \neq 0$, and
% $L_{1,P}$ be the projective line passing through $P$ and $\Omega_1$. Then
% $L_{1,P}$ is given by
% $L_{1,P}: Z_0Y - Y_0Z = 0$.

% Let $L_{2,P}$ be the line through $P$ and $\Omega_2$. Then $L_{2,P}$ is given by
% $L_{2,P}: Z_0X - X_0Z = 0$.
% \end{lemma}
% \begin{proof}
%   The points at infinity force $c_X=0$ and $c_Y=0$, respectively. The
%   equation then  follows from $P$ being a point on the line.\qed \end{proof}

%In the following we describe a special conic which passes through both points at
%infinity, $\Omega_1$ and $\Omega_2$, the point $\Ocal'$, and two arbitrary affine
%points $P_1$ and $P_2$ on $\EE{a}{d}$. 
Let $\phi(X,Y,Z) = c_{X^2}X^2 + c_{Y^2}Y^2 + c_{Z^2}Z^2 + c_{XY}XY + c_{XZ}XZ +
c_{YZ}YZ \in K[X,Y,Z]$ be a homogeneous polynomial of degree $2$ and 
$C: \phi(X,Y,Z) = 0$,
the associated plane (possibly degenerate) conic.
%
%\begin{remark}
  Since the points $\Omega_1,\Omega_2,\Ocal'$ are not on a line, a
 conic $C$ passing through these points cannot be a double line and  $\phi$ represents $C$
  uniquely up to multiplication by a scalar.
%\end{remark}
Evaluating $\phi$ at $\Omega_1, \Omega_2$, and $\Ocal'$, we see that a
conic $C$ through these points has the form
\begin{equation}\label{eqn:edpairings.edconic}
C: c_{Z^2}(Z^2 + YZ) + c_{XY}XY + c_{XZ}XZ = 0,
\end{equation}
where $(c_{Z^2}: c_{XY}: c_{XZ}) \in \proj^2(K)$.

\begin{theorem}\label{thm:edconic}
Let $\EE{a}{d}$ be a twisted Edwards curve over $K$, and let $P_1 = (X_1: Y_1 :
Z_1)$ and $P_2 = (X_2 : Y_2 : Z_2)$ be two affine, not necessarily distinct,
points on $\EE{a}{d}(K)$. Let $C$ be the conic passing through $\Omega_1$, $\Omega_2$,
$\Ocal'$, $P_1$, and $P_2$, \ie $C$ is given by an equation of the form~\eqref{eqn:edpairings.edconic}. If some of the
above points are equal, we consider $C$ and $\EE{a}{d}$ to
intersect with at least that multiplicity at the corresponding point. 
Then the
coefficients in \eqref{eqn:edpairings.edconic} of the equation $\phi$ of the conic $C$ are uniquely (up to scalars) determined as follows:
\begin{enumerate}[(a)]
\item
If $P_1 \neq P_2$, $P_1 \neq \Ocal'$ and $P_2 \neq \Ocal'$, then
\begin{eqnarray*}
c_{Z^2} & = & X_1 X_2 (Y_1 Z_2 - Y_2 Z_1),\\
c_{XY} & = & Z_1 Z_2 (X_1 Z_2 - X_2 Z_1 + X_1 Y_2 - X_2 Y_1),\\
c_{XZ} & = & X_2 Y_2 Z_1^2 - X_1 Y_1 Z_2^2 + Y_1 Y_2 (X_2 Z_1 - X_1 Z_2).
\end{eqnarray*}
\item
If $P_1 \neq P_2 = \Ocal'$, then
$c_{Z^2} = - X_1,\ c_{XY} = Z_1,\ c_{XZ} = Z_1$.
\item
If $P_1 = P_2$, then % \ne \Ocal'$,
\vspace*{-.6cm}
\begin{eqnarray*}
c_{Z^2} &=& X_1 Z_1 (Z_1-Y_1),\\
c_{XY} &=& d X_1^2Y_1 - Z_1^3, \\
c_{XZ} &=& Z_1(Z_1 Y_1 - aX_1^2).
\end{eqnarray*}
%T: this is included in (c)
% \item[(d)]If $P_1 = P_2=\Ocal'$,
% \begin{equation*}
% c_{Z^2}=0 ,\
% c_{XY}=1,\
% c_{XZ}=1.
% \end{equation*}
\end{enumerate}
\end{theorem}
\begin{proof}
%As for Lemma~\ref{lem:edpairings.conic1}, i
If the points are distinct, the
coefficients are obtained by evaluating the previous equation at the points
$P_1$ and $P_2$. We obtain two linear equations in $c_{Z^2},c_{XY}$, and $c_{XZ}$
\begin{eqnarray*}
c_{Z^2} (Z_1^2+Y_1Z_1)+c_{XY} X_1 Y_1 + c_{XZ} X_1Z_1 &=&0, \\
c_{Z^2} (Z_2^2+Y_2Z_2)+c_{XY} X_2 Y_2 + c_{XZ} X_2Z_2 &=&0. 
\end{eqnarray*}
The formulas in (a) follow from the (projective) solutions
\begin{equation*}
c_{Z^2} = \left|\begin{array}{cc} X_1 Y_1 & X_1 Z_1 \\ X_2 Y_2 & X_2 Z_2 \end{array}\right|,\, 
c_{XY} = \left|\begin{array}{cc} X_1 Z_1 & Z_1^2+Y_1 Z_1 \\ X_2 Z_2 &  Z_2^2+Y_2 Z_2 \end{array}\right|,\, 
c_{XZ} = \left|\begin{array}{cc} Z_1^2+Y_1 Z_1 & X_1 Y_1 \\ Z_2^2+Y_2 Z_2 & X_2 Y_2 \end{array}\right|. 
\end{equation*}

If $P_1=P_2 \ne \Ocal'$, we start by letting $Z_1=1,Z=1$ in the equations. The
tangent vectors at the non singular point $P_1=(X_1:Y_1:1)$ of $\EE{a}{d}$ and of $C$ are
$$\left(\begin{array}{c} d X_1^2 Y_1-Y_1 \\ a X_1-d X_1 Y_1^2 \end{array}\right), \quad
\left(\begin{array}{c} -c_{Z^2}-c_{XY} X_1 \\ c_{XY} Y_1+c_{XZ} \end{array}\right).$$
They are collinear if the determinant of their coordinates is zero which gives us a linear condition in the coefficients of $\phi$. We get a second condition by $\phi(X_1,Y_1,1)=0$. Solving the linear system, we get the projective solution
\begin{eqnarray*}
             c_{Z^2} &=& X_1^3 (-d Y_1^2+a)=X_1(1-Y_1^2)=X_1 (Y_1+1) (1-Y_1),\\
             c_{XY} &=& 2 d X_1^2 Y_1^2 -Y_1-Y_1^2+d X_1^2Y_1-a X_1^2\\
                    &=& -1 - Y_1 + dX_1^2Y_1^2 + dX_1^2Y_1 = (Y_1+1)(dX_1^2Y_1 - 1), \\
             c_{XZ} &=& -d X_1^2 Y_1^3-a X_1^2+Y_1^2+Y_1^3=(Y_1+1)(Y_1 - aX_1^2)
\end{eqnarray*}       
using the curve equation $a X_1^2+Y_1^2=1+d X_1^2 Y_1^2$ to simplify.
Finally, since $P_1 \ne \Ocal'$, we can divide by $1+Y_1$ and
homogenize to get the result which provides the formulas as stated.
% :
% \begin{eqnarray*}
% c_{Z^2} &=& X_1 Z_1 (Z_1-Y_1),\\
% c_{XY} &=& d X_1^2Y_1 - Z_1^3, \\
% c_{XZ} &=& Z_1(Z_1 Y_1 - aX_1^2).
% \end{eqnarray*}
%If $P_1 = P_2=\Ocal$ we may use the singular conic which is the
%product of the line $Y = Z$ tangent to $\EE{a}{d}$ in $\Ocal$ and the
%line $X = 0$ passing through the point $\Ocal'$. Thus $\phi= X(Z - Y)
%= XZ - XY$, so $c_{Z_2} = 0$, $c_{XY} = -1$, and $c_{XZ} = 1$, and the
%intersection multiplicity is at least $2$ in $\Ocal$. The 
%same values arise when evaluating the formulas under (c) at $P_1 =
%\Ocal$. Furthermore,
The same formulas hold if $P_1 = \Ocal'$ since
intersection multiplicity greater than or equal to $3$ at $\Ocal'$ is achieved by setting $\phi
= X(Y+Z) = XY + XZ$. 
%Not all three coefficients can be $0$,
%because this would imply $a = d$. 

Assume now that $P_1\neq P_2=\Ocal'$. Note that  the conic $C$ is tangent to
$\EE{a}{d}$ at $\Ocal'$ if and only if
$(\partial \phi/\partial x)(0,-1,1)=(c_{XY}y+c_{XZ}z)(0,-1,1)=0$, 
i.e.  $c_{XY}=c_{XZ}$. Then  $\phi= (Y+Z)(c_{Z^2} Z+c_{XY} X)$. Since $P_1 \ne
\Ocal'$, it is not on the line $Y+Z=0$. Then we get
$c_{Z^2} Z_1+c_{XY} X_1=0$ 
and the coefficients as in (b). 
\qed
\end{proof}

% T: the transformations come in the section on doubling; with all conditions.
% \begin{corollary}
% With the notation of Proposition~\ref{thm:edconic}, assume  that $P_1$, $P_2$ have either infinite order or order prime to $2$. Then we are either in case (a) or (c). Moreover we can replace the equations in case (c) by
% \begin{eqnarray*}
% c_{Z^2}&=&X_1(Y_1^2-Y_1Z_1),\\
% c_{XY}&=&Z_1(Z_1^2-aX_1^2-Y_1^2+Y_1Z_1),\\
% c_{XZ}&=&Y_1(aX_1^2-Y_1Z_1).
% \end{eqnarray*}
% \end{corollary}
% \begin{proof}
% For the last part of the statement. Since $P_1=P_2$ is not of order $4$, $Y_1 \ne 0$, so we can multiple the expressions in (c) by $-Y_1/Z_1$ and we get
% \begin{eqnarray*}
% c_{Z^2} &=& X_1Y_1(Y_1-Z_1), \\
% c_{XY} &=& -dX_1^2Y_1^2/Z_1 + Z_1^2Y_1 = -Z_1(aX_1^2 + Y_1^2 - Z_1^2) + Z_1^2Y_1,\\
% c_{XZ} &=& Y_1( aX_1^2-Z_1Y_1 ).
% \end{eqnarray*}
% \qed\end{proof}

%We now give the relation with the addition law. 
%The conic $C$
%described in Theorem~\ref{thm:edconic} gives a nice geometric
%interpretation of the group law on an Edwards curve, similar to the
%chord-and-tangent method of elliptic curves in Weierstrass form. \\

Let $P_1$ and $P_2$ be two affine $K$-rational points on a twisted Edwards curve
$\EE{a}{d}$,
and let $P_3 = (X_3:Y_3:Z_3) = P_1 + P_2$ be their sum. Let 
\begin{equation*}
l_1 = Z_3Y - Y_3Z,\quad l_2 = X
\end{equation*}
be the polynomials of the horizontal line $L_{1, P_3}$ through $P_3$ and the vertical line
$L_{2, \Ocal}$ through $\Ocal$ respectively, and let
\begin{equation*}
\phi = c_{Z^2}(Z^2 + YZ) + c_{XY}XY + c_{XZ}XZ 
\end{equation*}
be the unique polynomial (up to multiplication by a scalar) defined by
Theorem~\ref{thm:edconic}. %\\
The following theorem shows that the group law on a twisted Edwards curve indeed has a
geometric interpretation involving the above equations.  It gives
us an important ingredient to compute Miller functions.

\begin{theorem}\label{thm:edfunctions}
Let  $a,d \in K^*$ with $a \neq d$ and let
$\EE{a}{d}$ be a twisted Edwards curve over $K$. Let $P_1, P_2 \in \EE{a}{d}(K)$.
Define $P_3 = P_1 + P_2$. Let $\phi, l_1, l_2$ be defined as above. Then we have
\begin{equation} \label{edwardslaw}
\divi\left(\frac{\phi}{l_1l_2}\right) \sim (P_1) + (P_2) - (P_3) - (\Ocal).
\end{equation}
\end{theorem}
\begin{proof}
Let us consider the intersection divisor $(C \cdot \EE{a}{d})$ of the conic $C :
\phi=0$ and the singular quartic $\EE{a}{d}$.  Bezout's Theorem \cite[p. 112]{1969/fulton} tells us that the
intersection of $C$ and $\EE{a}{d}$ should have $2 \cdot 4=8$ points counting
multiplicities over $\overline{K}$. We note that the two points at infinity $\Omega_1$ and
$\Omega_2$ are singular points of multiplicity  $2$. Moreover, by definition of the conic $C$,   
$(P_1)+(P_2)+(\Ocal')+2 (\Omega_1)+2 (\Omega_2) \leq (C \cdot \EE{a}{d})$. Hence there is 
an eighth point $Q$ in the intersection. Let $L_{1,Q} : l_Q=0$ be the horizontal
line going through $Q$. Since the inverse for addition on twisted Edwards curves
is given by $(x,y) \mapsto (-x,y)$, we see that $(L_{1,Q} \cdot
\EE{a}{d})=(Q)+(-Q)-2(\Omega_2)$. On the other hand $(L_{2,\Ocal} \cdot
\EE{a}{d})=(\Ocal)+(\Ocal')-2(\Omega_1)$. Hence 
by combining the above divisors we get
$\divi\left(\frac{\phi}{l_Ql_2}\right) \sim (P_1) + (P_2) - (-Q) - (\Ocal)$.
By unicity of the group law with neutral element $\Ocal$ on the elliptic curve $\EE{a}{d}$
\cite[Prop.3.4]{1986/silverman}, the last equality means that $P_3=-Q$. Hence
$(L_{1,P_3} \cdot \EE{a}{d})=(P_3)+(-P_3)-2(\Omega_2)=(-Q)+(Q)-2(\Omega_2)$ and $l_1=l_Q$. So
$\divi\left(\frac{\phi}{l_1l_2}\right) \sim (P_1) + (P_2) - (P_3) - (\Ocal)$.
\qed
\end{proof}

\begin{remark}
From the proof, we see that $P_1+P_2$ is obtained as the 
mirror image with respect to the $y$-axis of the eighth intersection point 
of $\EE{a}{d}$ and the conic $C$ passing through $ \Omega_1,\Omega_2, \Ocal', P_1$ and $P_2$.
\end{remark}

\begin{example}\label{exmpl:edconic}
As an example we consider the Edwards curve $\EE{1}{-30}: x^2 + y^2 = 1 -
30x^2y^2$ over the set of real numbers $\R$. We choose the point
$P_1$ with $x$-coordinate $x_1 = -0.6$ and $P_2$ with $x$-coordinate $x_2 = 0.1$. 
Figure~\ref{fig:grouplawadd} shows addition of different
points $P_1$ and $P_2$, and Figure~\ref{fig:grouplawdbl} shows doubling of the
point $P_1$.
\begin{figure}[ht]
\centering
\psset{unit=1.5cm} 
\subfigure[$P_1 \neq P_2$, $P_1,P_2 \neq \Ocal'$, $P_3 = P_1 + P_2$\label{fig:grouplawadd}]{
\begin{pspicture}(-1.5,-1.5)(2.5,1.5) 
\psaxes[labels=none, ticks=none, linewidth=.2pt]{->}(0,0)(-1.5,-1.5)(1.5,1.5)  
\psplot[plotstyle=curve]{-1}{1}{1 x x mul sub 1 30 x mul x mul add div sqrt}
\psplot[plotstyle=curve]{-1}{1}{1 x x mul sub 1 30 x mul x mul add div sqrt neg}
\psdots(0,1)(0,-1)(-0.6,0.233)(0.1, 0.873)(-0.366306,0.415076)(0.366306,0.415076)
\psline[linestyle=dashed]{-}(-1.5,0.415076)(1.5,0.415076)
\uput[-55](-0.62,0.28){$P_1$}
\uput[45](0.1,0.873){$P_2$}
\uput[135](-0.366306,0.415076){$P_3$}
\uput[45](0.366306,0.415076){$-P_3$}
\uput[135](-1,0.415076){$L_{1,P_3}$}
\uput[-45](-1.4,0.26){$C$}
\uput[-45](0.4,-0.4){$\EE{1}{-30}$}
\uput[135](0,1){$\Ocal$}
\uput[-45](0,-1){$\Ocal'$}
\psplot[plotstyle=curve, linestyle=dashed, dash=10pt 1pt]{-1.5}{0.005}
{0.038386 0.369263 x mul add -0.038386 1.246886 x mul add div}
\psplot[plotstyle=curve, linestyle=dashed, dash=10pt 1pt]{0.065}{1.5}
{0.038386 0.369263 x mul add -0.038386 1.246886 x mul add div}
\end{pspicture}
}
\subfigure[$P_1 = P_2 \neq \Ocal'$, $P_3 = 2P_1$\label{fig:grouplawdbl}]{
\begin{pspicture}(-1.5,-1.5)(1.5,1.5) 
\psaxes[labels=none, ticks=none, linewidth=.2pt]{->}(0,0)(-1.5,-1.5)(1.5,1.5)  
\psplot[plotstyle=curve]{-1}{1}{1 x x mul sub 1 30 x mul x mul add div sqrt}
\psplot[plotstyle=curve]{-1}{1}{1 x x mul sub 1 30 x mul x mul add div sqrt neg}
\psdots(0,1)(0,-1)(-0.6,0.233)(-0.674653,-0.192817)(0.674653,-0.192817)
\psline[linestyle=dashed]{-}(-1.5,-0.192817)(1.5,-0.192817)
\uput[135](-0.6,0.233){$P_1$}
\uput[-105](-0.674653,-0.192817){$P_3$}
\uput[-45](0.674653,-0.192817){$-P_3$}
\uput[-135](-1,-0.192817){$L_{1,P_3}$}
\uput[45](-1.4,0.1){$C$}
\uput[45](0.4,0.4){$\EE{1}{-30}$}
\uput[45](0,1){$\Ocal$}
\uput[-45](0,-1){$\Ocal'$}
\psplot[plotstyle=curve, linestyle=dashed, dash=10pt 1pt]{-1.5}{-0.22}
{0.460267 0.127110 x mul add -0.460267 -3.515201 x mul add div}
\psplot[plotstyle=curve, linestyle=dashed, dash=10pt 1pt]{-0.04}{1.5}
{0.460267 0.127110 x mul add -0.460267 -3.515201 x mul add div}
\end{pspicture}
}
\caption{Geometric interpretation of the  group law on $x^2 + y^2 = 1 - 30x^2y^2$ over $\R$.\label{fig:edpairings.grouplaweq}}
\end{figure}
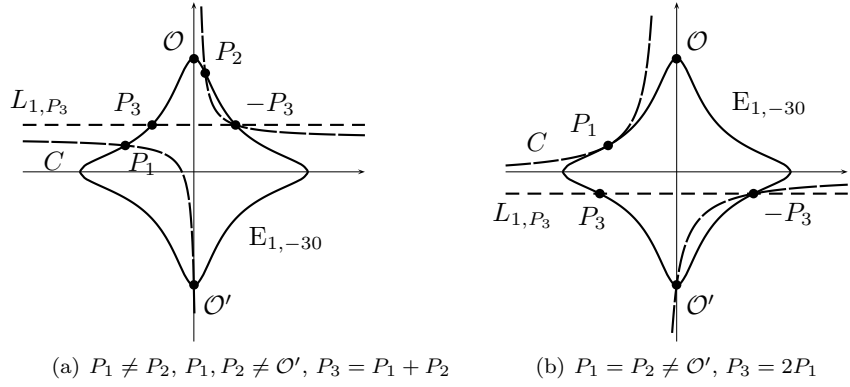
\end{example}

\vspace*{-1cm}
\section{Formulas for Pairings on Edwards Curves}\label{Miller}
In this section we show how to use the geometric interpretation of 
the group law 
to compute pairings. 
We assume that $k$ is even and that the 
second input point $Q$ is chosen by using the tricks in 
\cite{2002/barreto} and \cite{2004/barreto}:
%Note that, as explained in Section~\ref{edw}, on twisted Edwards curves $\EE{a}{d}$, twists affect the
%$x$-coordinate. 
Let $\F_{q^k}$ have basis $\{1,\alpha\}$ over
$\F_{q^{k/2}}$ with $\alpha^2=\delta \in \F_{q^{k/2}}$ and let $Q'=(X_0:Y_0:Z_0)\in
\EE{a\delta}{d\delta}(\F_{q^{k/2}})$.  Twisting $Q'$ with $\alpha$
ensures that the second argument of the pairing is on
$\EE{a}{d}(\F_{q^k})$ (and no smaller field) and is of the form
$Q=(X_0\alpha:Y_0:Z_0)$, where $X_0, Y_0, Z_0 \in \F_{q^{k/2}}$.

%% Check the loop shortening before submitting.
%This observation shows that our method can be used with loop shortening, i.e.
%to compute the twisted ate pairing.

By Theorem~\ref{thm:edfunctions} we have $g_{R,S}=\frac{\phi}{l_1l_2}$.
In each step of the Miller loop first
$g_{R,S}$ is computed,
it is then evaluated at $Q=(X_0\alpha:Y_0:Z_0)$ and finally $f$ is updated as $f\gets f
\cdot g_{R,P}(Q)$ (addition) or as $f\gets f^2\cdot g_{R,R}(Q)$ 
(doubling). Given the shape of
 $\phi$ and the point $Q=(X_0\alpha:Y_0:Z_0)$,
 we see that we need to compute
\begin{eqnarray*}
\frac{\phi}{l_1l_2}(X_0\alpha:Y_0:Z_0)
& = & \frac{c_{Z^2}(Z_0^2 + Y_0Z_0)+c_{XY}X_0\alpha Y_0+c_{XZ}X_0Z_0\alpha}
{(Z_3Y_0-Y_3Z_0)X_0\alpha}\\
& = & \frac{c_{Z^2}\frac{Z_0 + Y_0}{X_0\delta}\alpha+c_{XY}y_0+c_{XZ}}
{Z_3y_0-Y_3},\\
& \in  & (c_{Z^2}\eta\alpha+c_{XY}y_0+c_{XZ})\F^*_{p^{k/2}},
\end{eqnarray*}
\noindent
where $(X_3 : Y_3 : Z_3)$ are coordinates of the point $R+P$ or $R+R$, $y_0 =
Y_0/Z_0$, and $\eta=\frac{Z_0+Y_0}{X_0\delta}$.  Note that $\eta, y_0\in
\F_{q^{k/2}}$ and that they are fixed for the whole computation, so they
can be precomputed. 
The coefficients $c_{Z^2},c_{XY}$, and $c_{XZ}$ are defined over
$\F_q$, thus the evaluation at $Q$ given the coefficients of the conic
can be computed in $k \smults$ (multiplications by $\eta$
and $y_0$  need $\frac{k}{2} \smults$ each).

% In the next sections we give explicit formulas to efficiently compute
% $c_{Z^2}, c_{XY}$, and $c_{XZ}$ for addition and doubling.  For
% applications in cryptography we restrict our considerations to points
% in a group of prime order. Let the number of points on the curve
% factor as $\#\EE{a}{d}(\F_p)=4 hn$, with $n$ prime, 
% and let the basepoint $P$ have order
% $n$. This implies in particular that none of the additions or
% doublings involves $\Omega_1, \Omega_2$, or $\Ocal'$. The neutral
% element $\Ocal$ is a multiple of $P$, namely $nP$, but none of the 
% operations in the Miller loop will have it as its input. 
% This means that without
% loss of generality we can assume that none of the coordinates of the
% input points is $0$. In fact, for this assumption to hold we only need
% that $P$ has odd order, so that the points of order $2$ or $4$ are not
% multiples of it.

\subsection{Addition steps}\label{add}
Hisil et al. presented new addition formulas for twisted Edwards
curves in extended Edwards form 
 at Asiacrypt 2008~\cite{2008/hisil}. 
Let $P_3 =P_1 + P_2$ for two different points $P_1 = (X_1:Y_1:Z_1:T_1)$ and
$P_2 = (X_2:Y_2:Z_2:T_2)$ with $Z_1, Z_2\neq 0$ and $T_i=X_iY_i/Z_i$.
% for points $P=(X:Y:Z)$ with $Z\neq0$.  
%  in extended representation is given by
% \begin{eqnarray*}
% X_3 & = & (X_1Y_2 - Y_1X_2)(T_1Z_2 + Z_1T_2),\\
% Y_3 & = & (aX_1X_2 + Y_1Y_2)(T_1Z_2 - Z_1T_2),\\
% Z_3 & = & (aX_1X_2 + Y_1Y_2)(X_1Y_2 - Y_1X_2).
% \end{eqnarray*}
% We assume that the base
% point $P$ has odd order and so the condition $Z_1, Z_2\neq 0$ holds 
% during scalar 
% multiplication of $P$.
Theorem~\ref{thm:edconic} (a) %in Section~\ref{geom} 
states the coefficients
of the conic section for addition. We use $T_1, T_2$ to shorten the 
formulas.
\begin{eqnarray*}
c_{Z^2} & = & X_1 X_2 (Y_1 Z_2 - Y_2 Z_1)=Z_1Z_2(T_1X_2-X_1T_2),\\
c_{XY} & = & Z_1 Z_2 (X_1 Z_2 - Z_1X_2 + X_1 Y_2 - Y_1X_2),\\
c_{XZ} & = & X_2 Y_2 Z_1^2 - X_1 Y_1 Z_2^2 + Y_1 Y_2 (X_2 Z_1 - X_1 Z_2)\\
 & = & Z_1Z_2(Z_1T_2 - T_1Z_2 + Y_1 T_2 - T_1 Y_2).
\end{eqnarray*}

Note that all coefficients are divisible by $Z_1Z_2\neq 0$ and so we
scale the coefficients. The explicit formulas for computing $P_3 = P_1
+ P_2$ and $(c_{Z^2}, c_{XY},c_{XZ})$ are given as follows:
\begin{eqnarray*}
A&=&X_1\cdot X_2;\ B=Y_1\cdot Y_2;\ C=Z_1\cdot T_2;\ D=T_1\cdot Z_2;\ E=D+C;\\
F&=&(X_1-Y_1)\cdot (X_2+Y_2)+B-A;\ G=B+aA;\ H=D-C;\ I=T_1\cdot T_2;\\
c_{Z^2}&=&(T_1-X_1)\cdot (T_2+X_2)-I+A;\ c_{XY}=X_1\cdot Z_2-X_2\cdot Z_1+F;\\ 
c_{XZ}&=&(Y_1-T_1)\cdot (Y_2+T_2)-B+I-H;\\
X_3&=&E\cdot F;\ Y_3=G\cdot H;\ T_3=E\cdot H;\ Z_3=F\cdot G.
\end{eqnarray*}
With these formulas $P_3$ and $(c_{Z^2}, c_{XY},c_{XZ})$ can be
computed in $1\mults +(k+14)\smults + 1\smultsbya$, where $\smultsbya$ denotes the
costs of a multiplication by 
%the constant
$a$.  If the base point
$P_2$ has $Z_2 = 1$, the above costs reduce to $1\mults +(k+12)\smults +
1\smultsbya$.  We used Sage~\cite{2008/stein-sage} to verify the
explicit formulas. 
%  Note that there is no extra speed up from choosing
% $a=-1$ (unlike in \cite{2008/hisil}) since all subexpressions are also
% used in the computation of $(c_{Z^2}, c_{XY},c_{XZ})$.

\subsection{Doubling steps}\label{dbl}

Theorem~\ref{thm:edconic} (c) %in Section~\ref{geom}
 states the
coefficients of the conic section in the case of a doubling step. To speed up
the computation we multiply each coefficient by $-2Y_1/Z_1$;
remember that $\phi$ is unique up to scaling. Note also that $Y_1,
Z_1\neq 0$ because we assume that all points have odd order. The multiplication by
$Y_1/Z_1$ reduces the overall degree of the equations since we can use
the curve equation to simplify the formula for $c_{XY}$; the factor
$2$ is useful in obtaining an $\ssquarings-\smults$ tradeoff in the
explicit formulas below. We obtain:
\begin{eqnarray*}
c_{Z^2}&=&X_1(2Y_1^2-2Y_1Z_1),\\
c_{XY}&=&2(Y_1Z_1^3-dX_1^2Y_1^2)/Z_1=2(Y_1Z_1^3-Z_1^2(aX_1^2+Y_1^2)+Z_1^4)/Z_1\\
      &=&Z_1(2(Z_1^2-aX_1^2-Y_1^2)+2Y_1Z_1),\\
c_{XZ}&=&Y_1(2aX_1^2-2Y_1Z_1).
\end{eqnarray*}

Of course we also need to compute $P_3=2P_1$. We use the explicit
formulas from \cite{2008/BernsteinBirkner} for the doubling and reuse
subexpressions in computing the coefficients of the conic. The
formulas were checked for correctness with Sage~\cite{2008/stein-sage}.
%
% \begin{eqnarray*}
%   A&=&X_1^2;\ B=Y_1^2;\ C=Z_1^2; D=(X_1+Y_1)^2;\  E=(Y_1+Z_1)^2;\\
%   F&=&D-(A+B);\ G=E-(B+C);\ H=aA;\ I=H+B;\\
%   J&=&C-I;\ K=J+C;\ c_{XZ}=Y_1\cdot (2H-G);\ c_{XY}=Z_1\cdot (2J+G);\\
%   c_{Z^2}&=&F\cdot (Y_1-Z_1);\ X_3=F\cdot K;\  Y_3=I\cdot (B-H);\  Z_3=I\cdot K.
% \end{eqnarray*}
%
% These formulas compute $P_3$ %=(X_3:Y_3:Z_3)$ 
% and $(c_{Z^2},c_{XY},c_{XZ})$ in $6\smults +5\ssquarings + 1\smultsbya$. 
% If the doubling is followed by an addition the additional coordinate
% $T_3=X_3Y_3/Z_3$ needs to be computed. This is done by additionally
% computing $T_3=F\cdot(B-H)$ in $1 \smults$.
%
Since the input is given in extended form as $P_1=(X_1:Y_1:Z_1:T_1)$ we
can use $T_1$ in the computation of the conic as
\begin{eqnarray*}
c_{Z^2}&=&X_1(2Y_1^2-2Y_1Z_1)=2Z_1Y_1(T_1-X_1)
,\\
c_{XY}&=&Z_1(2(Z_1^2-aX_1^2-Y_1^2)+2Y_1Z_1),\\
c_{XZ}&=&Y_1(2aX_1^2-2Y_1Z_1)=2Z_1(aX_1T_1-Y_1^2),
\end{eqnarray*}
and then scale the coefficients by $1/Z_1$. The computation of
$P_3=(X_3:Y_3:Z_3:T_3)$ and $(c_{Z^2},c_{XY},c_{XZ})$ is then done in
$1\mults +1\squarings +(k+6)\smults +5\ssquarings + 2\smultsbya$ as
\begin{eqnarray*}
  A&=&X_1^2;\ B=Y_1^2;\ C=Z_1^2; D=(X_1+Y_1)^2;\  E=(Y_1+Z_1)^2;\\ 
  F&=&D-(A+B);\
  G=E-(B+C);\ H=aA;\ I=H+B;\ J=C-I;\\
  K&=&J+C;\ 
  c_{Z^2}=2Y_1\cdot (T_1-X_1);\ c_{XY}=2J+G;\ c_{XZ}=2(aX_1\cdot T_1-B);\\ 
  X_3&=&F\cdot K;\  Y_3=I\cdot (B-H);\  Z_3=I\cdot K;\ T_3=F\cdot(B-H).
\end{eqnarray*}
%
% These formulas compute $P_3=(X_3:Y_3:Z_3)$ and
% $(c_{Z^2},c_{XY},c_{XZ})$ in $6\smults +5\ssquarings + 2\smultsbya$.
% For computing the Tate pairing 
% this means
% that a doubling step costs $1\mults + 1\squarings + (k+6)\smults
% +5\ssquarings + 1\smultsbya$ in twisted Edwards coordinates and 
% $1\mults + 1\squarings + (k+6)\smults +5\ssquarings + 2\smultsbya$ 
% in extended coordinates.
%put something on the ate pairing here.

Note that like in \cite{2008/hisil} we can save $1\smultsbya$ per
doubling by changing to the extended representation only before an
addition.% Morain~\cite{2009/morain} showed that $a$ can
%always be chosen as $a=1$ when constructing a pairing-friendly
%curve. So the effective cost of multiplying by $a$ is $\smultsbya=0$.

\def\smultsbyc{{\bf m_{a_4}}} %multiplications by parameter d
\def\smultsbyb{{\bf m_{b}}} %multiplications by parameter b
\section{Operation counts}\label{comp}
We give an overview of the best formulas in the literature for computing the
Tate pairing
on Edwards curves and on the different forms of
Weierstrass curves in Jacobian coordinates. 
We compare the results with our new pairing formulas for Weierstrass
and Edwards curves.

Throughout this section we assume that $k$ is even, that 
the second input point 
$Q$ is given in affine coordinates, and that quadratic
twists are used so that multiplications with $\eta$ and $y_Q$ take 
$(k/2)\smults$ each.
\subsection{Overview}
Chatterjee, Sarkar, and Barua~\cite{2004/chatterjee} study pairings on
Weierstrass curves in Jacobian coordinates.  Their paper does not
distinguish between multiplications in $\F_q$ and in $\F_{q^k}$ but
their results are easily translated. For mixed addition steps their formulas
need $1\mults + (k+9)\smults +3\ssquarings$, and for doubling steps they need
$1\mults +(k+7)\smults+1\squarings+4\ssquarings$ if $a_4 = -3$.  For
doubling steps on general Weierstrass curves (no condition on $a_4$) the
formulas by Ionica and Joux~\cite{2008/ionica} are fastest with 
$1\mults +(k+1)\smults+1\squarings+11\ssquarings$.

Actually, any mixed addition step (mADD) or addition step (ADD) in Miller's
algorithm needs $1\mults
+k\smults$ for the evaluation at $Q$ and the update of $f$; each
doubling step (DBL) needs $1\mults +k\smults+1\squarings$ for the
evaluation at $Q$ and the update of $f$. In the following we do not
comment on these costs since they do not depend on the chosen
representation and are a fixed offset. We also do not report these
expenses in the overview table.

Hankerson, Menezes, and Scott~\cite{2009/hankerson} study pairing
computation on Barreto-Naehrig~\cite{2006/barreto} curves.  All BN
curves have the form $y^2=x^3+a_6$ and are thus more special than curves
with $a_4=-3$ or Edwards curves. 
% In their presentation they combine
%the pairing computation with the extension-field arithmetic and thus
%%the operations for the pure pairing computation are not stated
%explicitly but the formulas match those in \cite{2005/cheng}.  
They
need $6\smults +5\ssquarings$ for a doubling step and $9\smults +
3\ssquarings$ for a mixed addition step.
%  when computing the update functions
% for the Tate pairing.
Very recently, Costello et al. \cite{2009/costello} presented explicit formulas
for pairings
on curves of the form $y^2=x^3+b^2$, i.e. $a_4=0$ and $a_6$ is a square. 
Their representation is in projective rather than Jacobian coordinates.

To the best of our knowledge our paper is the first to publish full
(non-mixed) addition formulas for Weierstrass curves. Note that
\cite{2009/costello} started after our results became public.

%  compare our results with formulas in the literature, in particular
% with the pairing formulas for Edwards curves due to Ionica and
% Joux~\cite{2008/ionica} and the formulas by

%Note that the operation counts stated in \cite{2003/izu} are actually
%a bit worse than reported here. The speedups for even $k$ and
%denominator elimination were only introduced by Barreto et al.
%\cite{2002/barreto} in 2002 and so most older papers list $4k\smults$
%instead of $k\smults$ in the evaluation of the line function.  Often
%these papers go into details of doing arithmetic in the extension
%field $\F_{p^k}$ for very special $k$'s and state all operations in
%terms of operations in the base field.  We were unable to verify the
%operation counts in Table~1 of Cheng and Nistazakis~\cite{2005/cheng}
%for mixed addition.  Their formulas are identical to those in
%\cite{2003/izu} who give the correct operation counts; it seems that a
%mistake happened in copying data into the table.

% Our new formulas for Weierstrass curves in
% Section~\ref{appendix} are stated in the table as ``this paper''.

Das and Sarkar~\cite{2008/das} were the first to publish pairing
formulas for Edwards curves. We do not include them in our overview
since their study is specific to supersingular curves with $k=2$.
Ionica and Joux~\cite{2008/ionica} proposed the thus far fastest
pairing formulas for Edwards curves. Note that they actually compute
the $4$th power $T_n(P,Q)^4$ of the Tate pairing. This has almost no
negative effect for usage in protocols. So we include their result as
pairings on Edwards curves.

We denote Edwards coordinates by $\Ecal$, projective coordinates by
$\Pcal$, and Jacobian coordinates by $\Jcal$.
% The row ``$\Ecal$, this paper'' reports the results of the
% previous section using $2\Ecal\to \Ecal$ for the main doublings,
% $2\Ecal \to \Ecal^e$ for the final doubling, and $\Ecal^e+\Ecal^e\to
% \Ecal$ for the addition. 
Morain~\cite{2009/morain} showed that
 $2$-isogenies reach $a=1$ from any twisted Edwards curve; we therefore
 omit $\smultsbya$ in the table.\\

\begin{center}
\begin{tabular}{|l|l|l|l|}
\hline
&DBL&mADD&ADD\\\hline\hline
$\Jcal$, \cite{2008/ionica}, \cite{2004/chatterjee}&
$1\smults+11\ssquarings +1\smultsbyc$&$9\smults+3\ssquarings$&\hspace*{.3cm}---
\\\hline
$\Jcal$, \cite{2008/ionica}, this paper&
$1\smults+11\ssquarings +1\smultsbyc$&$6\smults+6\ssquarings$&$9\smults+6\ssquarings $
\\\hline
$\Jcal, a_4=-3$,
  \cite{2004/chatterjee}&
$7\smults+4\ssquarings$&$9\smults+3\ssquarings$&\hspace*{.3cm}---
\\\hline
$\Jcal, a_4=-3$,
this paper&
$6\smults+5\ssquarings$&$6\smults+6\ssquarings$&$9\smults+6\ssquarings$
\\\hline
$\Jcal, a_4=0$,
  \cite{2005/cheng}, \cite{2004/chatterjee}&
$6\smults+5\ssquarings$&$9\smults+3\ssquarings$&\hspace*{.3cm}---
\\\hline
$\Jcal, a_4=0$,
this paper&
$3\smults+8\ssquarings$&$6\smults+6\ssquarings$&$9\smults+6\ssquarings$
\\\hline
$\Pcal, a_4=0, a_6=b^2$ \cite{2009/costello}&
$3\smults+5\ssquarings$&$10\smults+2\ssquarings +1\smultsbyb$&
$13\smults+2\ssquarings +1\smultsbyb$
\\\hline
$\Ecal$, \cite{2008/ionica}&
$8\smults+4\ssquarings+1\smultsbyd$&$14\smults+4\ssquarings+1\smultsbyd$&\hspace*{.3cm}---
\\\hline
$\Ecal$, this paper&
$6\smults+5\ssquarings$&$12\smults$&$14\smults$
\\ \hline
\end{tabular}
\end{center}

\subsection{Comparison}
%It is common in the literature to assume $\ssquarings =0.8 \smults$
%but the $\ssquarings/\smults$ ratio actually depends on the shape of
%the prime $p$. For extremely sparse primes the reduction modulo $p$ is
%virtually free and so the difference in the speed of integer squarings
%and integer multiplication carries through to a difference in speed of
%finite field squarings and multiplications. If the prime is less
%sparse the reduction cost has more influence and the ratio will be
%closer to $1$.  The construction of pairing-friendly curves does not
%allow much freedom in choosing the prime $p$, so it rarely leads to
%very sparse primes $p$. To capture all situations we state the costs
%of the different formulas in terms of $\ssquarings$ and $\smults$
%separately.
%
% Note also that additions appear far less frequently
%than doublings since often the constructions lead to $n$ with low
%Hamming weight.

The overview shows that our new formulas for Edwards curves solidly
beat all previous formulas published for Tate pairing computation on
Edwards curves.  

Our new formulas for pairings on arbitrary Edwards curves are faster
than all formulas previously known for Weierstrass curves except for
the very special curves with $a_4=0$.  Specifically mixed additions on
Edwards curves are slower by some $\ssquarings-\smults$ tradeoffs but
doublings are much more frequent and gain at least an
$\ssquarings-\smults$ tradeoff each.  

The curves considered in \cite{2009/costello} are extremely special:
For $p\equiv 2 \bmod 3$ these curves are supersingular and thus have
$k=2$. For $p\equiv 1 \bmod 3$ a total of $3$ isomorphism classes
is covered by this curve shape. They have
faster doublings but slower additions and mixed additions than Edwards
curves.

Our own improvements to the doubling and addition formulas for
Weierstrass curves beat our new formulas for Edwards curves with
affine base point by several $\ssquarings-\smults$ tradeoffs. However,
in many protocols the pairing input $P$ is the output of some scalar
multiplication and is thus naturally provided in non-affine form.
Whenever converting $P$ to affine form is more expensive than proceeding in
non-affine form, all additions are full additions.  A full
addition on an Edwards curve needs one field operation less than on
Weierstrass curves. Depending on the frequency of addition and the
$\ssquarings/\smults$ ratio the special curves with $a_4=0$ might or
might not be faster. For all other curves, Edwards form is the
best representation. Furthermore, scalar multiplications on Edwards
curves are significantly faster than on Weierstrass curves.

Our new formulas for mixed addition steps (mADD) and doubling steps
(DBL) on Weierstrass curves are faster than all previous ones by
several $\ssquarings-\smults$ tradeoffs. Our formulas for full
addition (ADD) are the only ones in the literature for most
Weierstrass curves; for those with $a_4=0$ and $a_6=b^2$ they are
faster than those in \cite{2009/costello} for any
$\ssquarings/\smults$ ratio.

We note here that for curves in Weierstrass form the ate
pairing is more efficient than the Tate pairing, in particular when the R-ate
pairing or optimal pairings with a very short loop in Miller's algorithm are
computed, and when twists of degree $4$ and $6$ are used to represent torsion
points. Our comparison only refers to Tate pairing computation.

Further research needs to focus on how to compute variants of
the ate pairing on Edwards curves. To obtain the same or better efficiency as
the fastest pairings on Weierstrass curves, it needs to be clarified whether
optimal ate pairings can be computed and whether the above mentioned high-degree
twists can be used as well for suitable pairing-friendly curves in Edwards form.
Some initial results are presented in \cite{2010/CLN}.

\section{Construction of Pairing-Friendly Edwards Curves}\label{examples}
The previous chapter showed that pairing computation can benefit from
Edwards curves. Most constructions of pairing-friendly elliptic curves in the literature aim at a prime
group order and thus in particular do not lead to curves with
cofactor 4 that can be transformed to Edwards curves.  Galbraith,
McKee, and Valen\c{c}a \cite{2007/galbraith-mckee-valenca} showed how
to use the MNT construction to produce curves with small cofactor.
%These are the only curves with rho=1 that have group orders divisible
%by 4. For larger rho values it is more likely to construct Edwards
%curves.
%% This is not correct. The rho-value is already less than 1, if we have
%% cofactor 4.
Some other constructions that allow to find curves with cofactor divisible by
$4$ are described by Freeman, Scott, and Teske \cite{2008/taxonomy}.

To ensure security of the pairing based system two criteria must be
satisfied: The group $E(\F_p)$ must have a large enough prime order
subgroup so that generic attacks are excluded {\em and} $p^k$ must be
large enough so that index calculus attacks in $\F_{p^k}^*$ are
excluded.  For efficient implementation, we try to minimize $p$ and
$k$ to minimize the cost of arithmetic in $\F_p$ and $\F_{p^k}$ and
minimize $n$ to minimize the length of the Miller loop. This has the
effect of balancing the difficulty of the DLPs on the curve and in the
multiplicative group of the finite field $\F_{p^k}$.

Following the ECRYPT recommendations \cite{2009/ecrypt}, the ``optimal''
bitsizes of the primes $p$ and $n$ for curves $E/\F_p$ with $n \mid \#E(\F_p)$
and $n$ prime are shown in Table~\ref{tbl:bitsizes} for the most common security
levels. For these parameters, the DLP in the subgroup of
$E(\F_p)$ of order $n$ is considered equally hard as the DLP in $\F_{p^k}^*$.
In order to transform the curve to an Edwards curve, we 
need to have $\#E(\F_p) = 4hn$ for
some cofactor $h$. It follows that the rho-value $\rho = \log(p)/\log(n)$ of
$E$ is always larger than $1$. The recommendations imply a desired value for
$\rho\cdot k$ as displayed in Table~\ref{tbl:bitsizes}, which should be
achieved with an even embedding degree to favor efficient implementation.
This means that $p$ cannot be kept minimal but we managed to 
minimize $n$ to keep the Miller loop short.

In the following section we present six examples of pairing-friendly Edwards
curves with embedding degrees $k\in \{6,8,10,22\}$, which cover the
security levels given in Table~\ref{tbl:bitsizes}.  
%In Table~\ref{tbl:bitsizes} we assume $h=1$ and choose closest integers for $k$.
%Recall the definition of the rho-value $\rho = \log(p)/\log(n)$ of $E$.
%Note that in implementations we need even $k$ but for larger rho values 
%$\rho \cdot k$ plays the role of $k$.
% Using curves with a cofactor of the given size ensures that the prime $n$ is as
% small as possible for the corresponding size of $p$. The bit length of $n$
% is equal to the length of the Miller loop. Hence it can be minimized for
% the given security level by choosing a curve with parameter sizes as indicated in
% Table~\ref{tbl:bitsizes}.
%
\begin{table}[h!]
\centering\begin{tabular*}{8cm}{@{\extracolsep{\fill}}|r||r|r|r|r|r|r|}
\hline
security & $80$ & $96$ & $112$ & $128$ & $160$ & $256$ \\
\hline\hline
$\log_2(n)$ & $160$ & $192$ & $224$ & $256$ & $320$ & $512$ \\
\hline
$\log_2(p^k)$ & $1248$ & $1776$ & $2432$ & $3248$ & $4800$ & $15424$ \\
\hline
$\rho\cdot k$ &  $7.80$ &    $9.25$ & $10.86$ & $12.67$ &  $15$ & $30.13$\\
\hline
\end{tabular*}
\vspace{0.2cm}
\caption{``Optimal'' bitsizes for the primes $n$ and $p$ and the corresponding
values for $\rho\cdot k$ for most common security levels.\label{tbl:bitsizes}}
\end{table}
\vspace{-1cm}
%
%
% \begin{table}[h!]
% \centering\begin{tabular*}{5cm}{@{\extracolsep{\fill}}|r|r|r|}
% \hline
% $n$ & $p^k$ & $\rho\cdot k$\\
% \hline\hline
% $160$  & $1248$ & $7.80$\\
% \hline
% $192$  & $1776$ & $9.25$\\
% \hline
% $224$  & $2432$ & $10.86$\\
% \hline
% $256$  & $3248$ & $12.67$\\
% \hline
% $320$ & $4800$ & $15$ \\
% \hline
% $512$  & $15424$ & $30.13$\\
% \hline
% \end{tabular*}
% %\vspace{0.2cm}
% \caption{``Optimal'' bitsizes for the primes $n$ and $p$ and the corresponding
% values for $\rho\cdot k$. \label{tbl:bitsizes} }
% \end{table}
% %
%

%\bibliography{bib}
%\bibliographystyle{plain}

\section{Examples of Pairing-Friendly Edwards Curves}
This section presents pairing-friendly Edwards curves. Note that they 
were constructed for applications using the Tate pairing so that the curve
over the ground field has a point of order 4. %(rather than for the ate 
%pairing for which the twist of the curve should be birationally equivalent 
%to a curve in Edwards form.
They are all defined over a prime field $\F_p$, and the $\rho$ values are stated with the
curves.  Notation is as before, where the number of $\F_p$-rational
points on the curve is $4hn$.

The curve examples in this section cover the security levels in
Table~\ref{tbl:bitsizes}. We used the method and formula in \cite{2009/ecrypt}
to determine the effective security in bits on the curve and in the finite
field.

\parindent=0em
%\begin{itemize}
% \item $D=1$, $\lceil\log(n)\rceil = 363$, $\lceil\log(h)\rceil = 7$, $\lceil\log(p)\rceil = 371$
% {\small
% \begin{eqnarray*}
% p & = & 32428903728427434871960638456028409162281939582432575945\\
%   &   & 30632153559402628010019946681624958973937239637420169141,\\
% n & = & 11105788948091587284918026868502879850096554651518005460\\
%   &   & 623832064312035897815509951488907964532000965993787241,\\
% h & = & 73,\\
% d & = & 16214451864213717435980319228014204581140969791216287972\\
%   &   & 65316076779701314005009973340812479486968619818710084571.
% \end{eqnarray*}
% }

%\item 
\subsubsection{Security level 80 bits (generic: 82 bits, index calculus: 79 bits):\\}
$k=6, \rho=1.22$ following \cite{2007/galbraith-mckee-valenca}:\vspace{0.2cm}\\
  $D = 7230$, $\lceil\log(n)\rceil=165, \lceil\log(h)\rceil=34,
  \lceil\log(p)\rceil=201$, $k\lceil \log(p) \rceil = 1206$%\\
{\small
\begin{eqnarray*}
p & = & 2051613663768129606093583432875887398415301962227490187508801,\\
n & = & 44812545413308579913957438201331385434743442366277,\\
h & = & 7\cdot733\cdot2230663,\\
d & = & 1100661309421493056836745159318889208210931380459417578976626.%\\
\end{eqnarray*}
}

%\item 
\subsubsection{Security level 96 bits (generic: 95 bits, index calculus: 93 bits):\\}
$k=6, \rho=1.48$ following \cite{2007/galbraith-mckee-valenca}:\vspace{0.2cm}\\
 %\\
  $D = 4630$, $\lceil\log(n)\rceil=191, \lceil\log(h)\rceil=90,
  \lceil\log(p)\rceil=283$, $k\lceil \log(p) \rceil = 1698$%\\
{\small
\begin{eqnarray*}
p & = & 1207642247325762099962277292422023053565510428560082635785
% \\
%   &   & 
6070179619031510615886361601,\\
n & = & 2498886235887409414948289020220476887707263210939845485839,\\
h & = & 11161\cdot19068349\cdot5676957216676051,\\
d & = & 2763915426899189358845059350727381504946815286189972438681082636399984067165911590884.
\end{eqnarray*}
}

%\item $k=10, \rho=1.49$ following Construction~6.5 in \cite{2008/taxonomy}:
%  $D = 1$, $\lceil\log(n)\rceil=328, \lceil\log(h)\rceil=160, \lceil\log(p)\rceil=490$\\
%%using x=2199023267236
%  {\small
%\begin{eqnarray*}
%p &=& 985055064842915992056265682718114432717943699383981706103022332940240876 \\
% &  &4103875127551514361113305839062859255153601,\\ 
%n & =& 115792474010234409409945101324928128192675797160629543631596278493065052913601 \\
%h &=& 21267683268341151449717030304934560000 \\
%d &=& -1.
%\end{eqnarray*}}

%\item
%\newpage
\subsubsection{Security level 112 bits (generic: 112 bits, index calculus: 117 bits):\\}
$k=8, \rho=1.50$ following Example~6.10 in \cite{2008/taxonomy}:\vspace{0.2cm}\\
  $D = 1$, $\lceil\log(n)\rceil= 224, \lceil\log(h)\rceil= 111,
  \lceil\log(p)\rceil= 337$, $k\lceil \log(p) \rceil = 2696$%\\
  {\small
\begin{eqnarray*}
p & = & 
2337736653699105669260383900156918881424547469292956866896259132890909437035723\\
 & & 48756028778874481604289\\
n & = & 22985796260053765810955211899935144604417092746113717429138553265289\\
h & = & 315669989 \cdot 558193107149 \cdot 14429732414341\\
d & = &
2137384144163601288355195724634322855348958454823252387999763620028079615999998\\
& & 48556640836158104712032\\
\end{eqnarray*}
}

\subsubsection{Security level 128 bits (generic: 133 bits, index calculus: 127 bits):\\}
$k=8, \rho=1.50$ following Example~6.10 in \cite{2008/taxonomy}:\vspace{0.2cm}\\
  $D = 1$, $\lceil\log(n)\rceil=267, \lceil\log(h)\rceil=133,
  \lceil\log(p)\rceil=401$, $k\lceil \log(p) \rceil = 3208$%\\
  {\small
\begin{eqnarray*}
p & = & 
5106500003052745062671102775396566649855857676935384847563820321458497449535443\\
 & & 6071209268470508469629312810691036880709,\\
n &= & 
8337030425086788445100704671763896482549397437850042633596560118040562641504433,\\
h &=& 
5 \cdot 17 \cdot 1229 \cdot 3181 \cdot 4608053164778689785613892277341,\\
d &=& 
2553250001526372531335551387698283324927928838467692423781910160729248724767721\\
 & & 8035604634235254234814656405345518440355,\\
\end{eqnarray*}
}

%\item 
\subsubsection{Security level 160 bits (generic: 164 bits, index calculus: 154 bits):\\}
$k=10, \rho=1.49$ following Construction~6.5 in
\cite{2008/taxonomy}:\vspace{0.2cm}\\
  $D = 1$, $\lceil\log(n)\rceil=328, \lceil\log(h)\rceil=160,
  \lceil\log(p)\rceil=490$, $k\lceil \log(p) \rceil = 4900$%\\
%using x=2199023267236
  {\small
\begin{eqnarray*}
p & = & 319667071934078971315677746964738362812713703914060344412320604868708613896665173327525\\
& &2543330209754427990875101879841425427646115157594515629491249,\\
n & = & 546812704438652190176048473638362779688423061794499756311925945545462152449512232744941\\
& &959488864241,\\
h & = & 2^4 \cdot 70199^4 \cdot 7831391^4,\\
d & = &
366838958032886838857360394166535857747556934852621175164120734346101628194129743602008\\
& & 259319768868802620569094456792293200142806009932471922115210.
\end{eqnarray*}
}

%\item 
\subsubsection{Security level 256 bits (generic: 259 bits, index calculus: 259 bits):\\}
$k=22, \rho=1.39$ following Construction~6.6 in
\cite{2008/taxonomy}:\vspace{0.2cm}\\
  $D = 3$, $\lceil\log(n)\rceil=519, \lceil\log(h)\rceil=204,
  \lceil\log(p)\rceil=724$, $k\lceil \log(p) \rceil =15928 $%\\
%using x=2199023267236
  {\small
\begin{eqnarray*}
p &=& 793243907836538225101919663581953770913765580662849594203574636874518836858270555160144\\
  & & 920983827280386815433912190214824741372960533715598691121880716182459140439367767771926\\
  & & 66177113943586415044911851669785290654695123,\\
n &=& 962131187808560377898569195262572710988984869464755002509459666178069262628367282191252\\
  & &973105101373704953818660670550658659790389637917606342501732923486369,\\
h &=& 3^5 \cdot 7 \cdot 13^2 \cdot 19^2 \cdot 37^2 \cdot 6421^2 \cdot 7219 \cdot 3498559^2 \cdot 22526869^2 \cdot 78478074679,\\
d &=& 264414627547939780810839826727395383259987444981352560753582877086320074680650633780571\\
  & &920373615518032509200852332864216413041328949865016666759728218019456097204687710831048\\
& &17656092016879614901160245443945786256399518.
\end{eqnarray*}
}

% \item $D = 314$, $\lceil\log(n)\rceil=220, \lceil\log(h)\rceil=98, \lceil\log(p)\rceil=319$\\
% {\small
% \begin{eqnarray*}
% p & = & 9452707311247707513330618188853923205411626343551115070653265144510\\
%   &   & 04408844212168675659778272001,\\
% n & = & 1336495861025991472146331033760710418580743090769112585053164534599,\\
% h & = & 3\cdot58939622055090151702905271933,\\
% d & = & 6987035440681854303189570169230960965572653735366548738901481656990\\
%   &   & 62950896708357057310497167360.\\
% \end{eqnarray*}
% }
%\end{itemize}

\bibliography{bib}
\bibliographystyle{plain}

\end{document}